\documentclass[12pt]{article}

\usepackage{graphicx}
\usepackage{amsmath}

\usepackage{mathptmx}      

\usepackage{bm}
\usepackage{amssymb}
\usepackage{amsthm}
\usepackage{tikz}
\usepackage{pgfplots}
\usepackage{verbatim}

\makeatletter
\renewcommand{\@biblabel}[1]{#1.} 
\makeatother

\begin{document}

\title{A priori estimation of a time step for numerically solving parabolic problems}

\author{Petr N. Vabishchevich\footnote{This work was supported by RFBR (project 13-01-00719)}}

\date{\small Nuclear Safety Institute, Russian Academy of Sciences, 52, B. Tulskaya, Moscow, Russia \\
      North-Eastern Federal University, 58, Belinskogo, Yakutsk, Russia \\
      E-mail: vabishchevich@gmail.com}

\maketitle

\begin{abstract}

This work deals with the problem of choosing a time step for the numerical solution
of boundary value problems for parabolic equations.
The problem solution is derived using the fully implicit scheme, whereas a time step is selected via
explicit calculations. The selection strategy consists of the following steps.
First, using the explicit scheme, we calculate the solution at a new time level.
Next, we employ this solution in order to obtain  the solution
at the previous time level (the implicit scheme, explicit calculations).
This solution should be close to the solution of our problem at this time level with a prescribed accuracy.
Such an algorithm leads to explicit formulas for the calculation of the time step
and takes into account both the dynamics of the problem solution and changes in coefficients of the equation
and in its right-hand side.
The same formulas for the evaluation of the time step we get using a comparison of two
approximate solutions, which are obtained using the  explicit scheme with 
the primary time step and the step that is reduced by half.
Numerical results are presented for a model parabolic boundary value problem,
which demonstrate the robustness of the developed algorithm for the time step selection.


\end{abstract}

\thispagestyle{empty}

\section{Introduction}

In numerically solving  boundary value problems for time-dependent equations, 
emphasis is on discretizations in time \cite{Angermann,Ascher,LeVeque}.
For parabolic equations of second order, unconditionally stable schemes are constructed using
implicit approximations \cite{Samarskii,SamarskiiGulin,SamarskiiMatusVabishchevich}.
Two-level schemes are commonly used in computational practice, whereas multilevel schemes
occur more rarely. For unconditionally stable schemes, a time step is 
selected only due to the accuracy of the approximate solution.

The problem of the control over a time step is relatively well resolved for
the numerically solving Cauchy problem for systems of differential equations
\cite{ascher1998computer,Gear1971,HairerNorsettWanner1987}.
The basic approach involves the following stages. First, we perform additional calculations
in order to estimate the error of the the approximate solution at a new time level.
Further, a time step is estimated using the theoretical asymptotic dependence of accuracy on a time step.
After that we decide is it necessary to correct the time step and to repeat calculations.

Additional calculations for estimating the error of the approximate solution
may be performed in a different way. In particular, it is possible to obtain an approximate solution
using two different schemes that have the same theoretical order of accuracy.
The most famous example of this strategy involves the solution of the problem on a separate time interval
using a preliminary step (the first solution ) and the step reduced by half (the second solution).
In numerically solving the Cauchy problem for systems of ordinary differential equations,
there are are also applied nested methods, where two approximate solutions of different
orders of accuracy are compared.

In the above-mentioned methods of selecting a time step, a posteriori estimation of
accuracy is employed. In this case, we decide is this time step acceptable or it is necessary to change it  for re-calculations (increase or reduced and how much) only after performing calculations.
Such strategies can be also applied to the approximate solution of unsteady boundary value problems
using a more advanced a posteriori analysis
\cite{bangerth2003adaptive,moller2011adaptive,verfurth2013posteriori}.

In this paper, we consider an a priori selection of a time step for the approximate
solution of boundary value problems for parabolic equations.
To obtain the solution at a new time level, the backward Euler scheme is employed.
The time step at the new time level is explicitly calculated using two
previous time levels and takes into account changes in the equation
coefficients and its right-hand side.
The paper is organized as follows.
In Section 2, we consider a Cauchy problem for a system of linear ordinary
differential equations that is obtained from numerically solving
boundary value problems for parabolic equations after discretization in space.
For the approximate solution, estimates for stability are presented along with
estimates for accuracy in the corresponding Hilbert space.
Formulas for the selection  of a time step are obtained in Section 3 using a
comparison of the problem solutions corresponding to the forward  time level and backward one.
In Section 4, we show that similar expressions for a time step can be obtained via
making a comparison of the solutions derived with one time step and two half steps.
Section 5 presents numerical results for a model boundary value problem for a one-dimensional parabolic equation obtained on the basis of the developed algorithm for selecting a time step.

\section{Model problem} 

Let us consider the Cauchy problem for the linear equation
\begin{equation}\label{2.1}
  \frac{du}{dt} + A(t) u = f(t),
  \quad 0 < t \leq T, 
\end{equation}
supplemented with the initial condition 
\begin{equation}\label{2.2}
  u(0) = u_0 .
\end{equation}
The problem is investigated in a finite-dimensional Hilbert space $H$.
Assume that
\[
 A(t) \geq 0 
\] 
in $H$.
Due to the non-negativity of the operator $A$, for the problem (\ref{2.1}), (\ref{2.2}),
we have the following estimate for stability with respect for the initial data and the right-hand side:
\begin{equation}\label{2.3}
  \|u(t)\| \leq \|u_0\| + \int_{0}^{t} \|f(\theta) \| d \theta .
\end{equation} 

The problem (\ref{2.1}), (\ref{2.2}) results from finite difference, finite volume or finite element approximations (lumped masses scheme \cite{thomee2010galerkin}) for numerically solving boundary value 
problems for a parabolic equation of second order. In this problem, an unknown function $u({\bm x},t)$
satisfies the equation
\[
   \frac{\partial u}{\partial t} 
   - \sum_{\alpha =1}^{m}
   \frac{\partial }{\partial x_\alpha} 
   \left ( k({\bm x},t)  \frac{\partial u}{\partial x_\alpha} \right ) + c({\bm x},t) u = f({\bm x},t),
   \quad {\bm x}\in \Omega,
   \quad 0 < t \leq  T,
\]
where 
$\underline{k} \leq k({\bm x}) \leq  \overline{k}, \ {\bm x} \in \Omega$,
$\underline{k} > 0$, $c({\bm x},t) \geq 0$.
The equation is complemented by the Dirichlet boundary conditions
\[
   u({\bm x},t) = g({\bm x},t),
   \quad {\bm x}\in \partial \Omega,
   \quad 0 < t \leq  T,
\]
and the initial condition
\[
   u({\bm x},0) = u_0({\bm x}),
   \quad {\bm x}\in \Omega.
\]

To solve numerically this time-dependent problem, we introduce a non-uniform grid in time:
\[
 t_0=0, \quad t_{n+1} = t_n + \tau_{n+1},
 \quad n = 0,1, ..., N-1,
 \quad t_N = T .   
\]
We will employ notation $f_n = f(t_n)$. For the problem (\ref{2.1}), (\ref{2.2}), we apply
the fully implicit scheme, where the transition from the current time level to the next one
is performed as follows:
\begin{equation}\label{2.4}
  \frac{y_{n+1} - y_{n}}{\tau_{n+1}} + A_{n+1} y_{n+1} = f_{n+1},
  \quad n = 0,1, ..., N-1, 
\end{equation} 
starting from the initial condition
\begin{equation}\label{2.5}
 y_0 = u_0.
\end{equation} 

Under the restriction  $A_{n+1} \geq 0$, from (\ref{2.4}), it follows immediately that
the approximate solution satisfies the level-wise estimate
\[
 \|y_{n+1}\| \leq \|y_{n}\| + \tau_{n+1} \|f_{n+1}\| .
\]
Thus, we obtain the discrete analog of the estimate (\ref{2.3}):
\begin{equation}\label{2.6}
 \|y_{n+1}\| \leq \|u_{0}\| + \sum_{k=0}^{n} \tau_{k+1} \|f_{k+1}\| 
\end{equation} 
corresponding to the problem (\ref{2.4}), (\ref{2.5}).
For the error  $z_n = y_n - u_n$ of the approximate solution, we have the problem
\[
  \frac{z_{n+1} - z_{n}}{\tau_{n+1}} + A_{n+1} z_{n+1} = \psi_{n+1},
  \quad n = 0,1, ..., N-1,  
\] 
\[
 z_0 = 0.
\] 
Here $\psi_{n+1}$ stands for the truncation error:
\begin{equation}\label{2.7}
 \psi_{n+1} = f_{n+1} -
 \frac{u_{n+1} - u_{n}}{\tau_{n+1}} + A_{n+1} u_{n+1} . 
\end{equation} 
Similarly to (\ref{2.6}), we get the estimate for error:
\[
  \|z_{n+1}\| \leq \sum_{k=0}^{n} \tau_{k+1} \|\psi_{k+1}\| .
\]
Therefore, to control the error, we can employ the summarized error
$\tau_{n+1} \delta$ over the interval  $t_n \leq t \leq t_{n+1}$.
In this case, a value $\delta$ defines the same level of the error over the entire interval of integration.

\section{Algorithm for estimation of a time step} 

If we will be able to calculate the truncation error $\psi_{n+1}$, then it will be possible to get a posteriori
estimate of the error. Comparing $\|\psi_{n+1}\|$ with the prescribed error level $\delta$,
this makes possible to evaluate the quality of the choice of the time step $\tau_{n+1}$.
Namely, if $\|\psi_{n+1}\|$ is much larger (smaller)  than $\delta$, then the time step is taken too large (small), and if $\|\psi_{n+1}\|$ is close to $\delta$, then this time step is optimal.
Thus, we have
\begin{equation}\label{3.1}
  \tau_{n+1}: \ \|\psi_{n+1}\| \approx \delta .
\end{equation}
The problem is that we cannot evaluate the truncation error, since it is determined
using the exact solution that is unknown. Because of this, we must focus
on some estimates for the truncation error that guarantee the fulfilment of (\ref{3.1}). 

The following strategy is proposed for the correction of the time step.
The step $\tau_{n+1}$ is selected from the conditions:
\begin{description}
 \item[Forward step.] Using the explicit scheme, we calculate the solution $v_{n+1}$ at the time level
  $t_{n+1}$;
 \item[Backward step.] From the obtained $v_{n+1}$, applying the implicit scheme, we determine $v_{n}$ at the     time level $t_{n}$ (explicit calculations);
 \item[Step selecting.] The step $\tau_{n+1}$ is evaluated via closeness between $v_{n}$ and $y_{n}$.
\end{description}
In fact, we carry out the back analysis of the error of the approximate solution over the interval
$t_n \leq t \leq t_{n+1}$ using two schemes (explicit and implicit) of the same accuracy. 

Let us present the formulas for selecting a time step.
The solution  $v_{n+1}$ is determined from the equation
\begin{equation}\label{3.2}
  \frac{v_{n+1} - y_{n}}{\tau_{n+1}} + A_{n} y_{n} = f_{n} .
\end{equation} 
For $v_{n}$, we have
\begin{equation}\label{3.3}
  \frac{v_{n+1} - v_{n}}{\tau_{n+1}} + A_{n+1} v_{n+1} = f_{n+1} .
\end{equation} 
From (\ref{3.2}), (\ref{3.3}), we immediately get
\begin{equation}\label{3.4}
  v_{n} - y_{n} = \tau_{n+1} (A_{n+1} - A_{n}) y_{n} - \tau_{n+1} (f_{n+1} - f_n) +
  \tau_{n+1}^2 A_{n+1} (f_n - A_ny_n) .
\end{equation} 

The first two terms are associated with the time derivative applied to the problem operator and to the right-hand side. To evaluate them approximately, it seems reasonable to use the time step from the previous time level. But this may be inconvenient to implement.

For instance, we have 
\[
 \tau_{n+1} (f_{n+1} - f_n) = \tau_{n+1}^2 \frac{f_{n+1} - f_n}{\tau_{n+1}},
\]
and therefore we have to evaluate the difference derivative of the right-hand side
for $t_{n} \leq t \leq t_{n+1}$.
The problem is that the derivation of such estimates involves the unknown value  $t_{n+1}$.
The simplest approach is to evaluate this derivative using the previous time step:
\[
 \frac{f_{n+1} - f_n}{\tau_{n+1}} \approx \frac{f(t_n + \tau_n) - f_n}{\tau_{n}} .
\]
But in this case, if $\tau_{n+1} > \tau_{n}$ , then we cannot detect significant changes
in the right-hand side for $t_{n} + \tau_{n} \leq t \leq t_{n} + \tau_{n+1}$.

To resolve the problem, it is possible to use the standard methods available to control a time step for numerically solving time-dependent problems.
The first method restricts the growth of the time step with respect to the previous value.
We set
\begin{equation}\label{3.5}
 \tau_{n+1} \leq \gamma \tau_{n}, 
\end{equation} 
where $\gamma > 1$ is a numerical parameter.
The second requirement is that the step cannot be too small:
\begin{equation}\label{3.6}
 \tau_{n+1} \geq  \tau_{0}, 
\end{equation}  
where $\tau_{0}$ is a specified minimum time step.

Under the assumption  (\ref{3.5}), we can estimate the time derivative of the right-hand side,
putting
\[
 \frac{f_{n+1} - f_n}{\tau_{n+1}} \approx \frac{f(t_n + \gamma \tau_n) - f_n}{\gamma \tau_{n}} .
\] 
Therefore
\[
 \tau_{n+1} (f_{n+1} - f_n) \approx \frac{\tau_{n+1}^2}{\gamma\tau_{n}} (\widetilde{f}_{n+1} - f_{n}) ,
\] 
\[
 \tau_{n+1} (A_{n+1} - A_{n}) y_{n} \approx  \frac{\tau_{n+1}^2}{\gamma\tau_{n}} (\widetilde{A}_{n+1} - A_{n}) y_{n}, 
\] 
where
\begin{equation}\label{3.7}
 \widetilde{f}_{n+1} = f(t_n + \gamma\tau_n),
 \quad  \widetilde{A}_{n+1} = A(t_n + \gamma\tau_n) . 
\end{equation}
For the last term in the right-hand side of (\ref{3.4}), in view of (\ref{2.4}), we have
\[
 \tau_{n+1}^2 A_{n+1} (f_n - A_ny_n) = \frac{\tau_{n+1}^2}{\tau_{n}}  A_{n+1} (y_{n}-y_{n-1}) .
\]
With accuracy up to $O(\tau_{n+1}^3)$, we put 
\[
 \tau_{n+1}^2 A_{n+1} (f_n - A_ny_n) \approx  \frac{\tau_{n+1}^2}{\tau_{n}}  \widetilde{A}_{n+1} (y_{n}-y_{n-1}) .
\]
With this in mind, the equality (\ref{3.4}) is replaced by the approximate equality:
\begin{equation}\label{3.8}
  v_{n} - y_{n} \approx \frac{\tau_{n+1}^2}{\tau_{n}} \left (  
  \frac{1}{\gamma} (\widetilde{A}_{n+1} - A_{n}) y_{n} - \frac{1}{\gamma} (\widetilde{f}_{n+1} - f_n) +
 \widetilde{A}_{n+1} (y_{n}-y_{n-1}) \right ) .
\end{equation} 

The value of $v_{n} - y_{n}$ we associate with the solution error over the interval
$t_n \leq t \leq t_{n+1}$. Because of this, we set
\begin{equation}\label{3.9}
  \|v_{n} - y_{n} \| \leq \tau_{n+1} \delta .
\end{equation} 
From (\ref{3.8}), we have
\begin{equation}\label{3.10}
 \|v_{n} - y_{n} \| \lesssim \frac{\tau_{n+1}^2}{\tau_{n}} \left (  
  \frac{1}{\gamma}\|\widetilde{A}_{n+1} - A_{n}) y_{n} \| + \frac{1}{\gamma}\|\widetilde{f}_{n+1} - f_n \| +
 \| \widetilde{A}_{n+1} (y_{n}-y_{n-1}) \| \right ) .
\end{equation}

In view of (\ref{3.5}), (\ref{3.6}), (\ref{3.9}),  from (\ref{3.8}), we obtain
the following formula for calculating the time step:
\begin{equation}\label{3.11}
\begin{split}
  \tau_{n+1} & = \max\left \{\tau_0, \min\{\gamma, \gamma_{n+1}\} \tau_{n} \right \} , \\
  \gamma_{n+1} & = \delta \left (  \frac{1}{\gamma}\| (\widetilde{A}_{n+1} - A_{n}) y_{n} \| 
  + \frac{1}{\gamma}\|\widetilde{f}_{n+1} - f_n \| +
 \| \widetilde{A}_{n+1} (y_{n}-y_{n-1}) \| \right )^{-1}  .
\end{split}
\end{equation}
This formula for selecting a time step reflects clearly (see the denominator in the expression for $\gamma_{n+1}$) corrective actions,  which are related to the time-dependence of the problem operator (the first part) and the right-hand side (the second part) as well as to the time-variation of the solution itself (the third part).
 
\section{Estimation of a time step on the basis of step doubling} 

To solve numerically the Cauchy problem, the traditional strategy is to select an integration step
using a comparison of the approximate solution obtained by the preliminary step with the solution calculated with the step reduced by half.
For numerically solving problem (\ref{2.1}), (\ref{2.2}), we use fully implicit scheme
(\ref{2.4}), (\ref{2.5}). We employ the explicit scheme over the interval $t_n \leq t \leq t_{n+1}$
in order  to select the time step $\tau_{n+1}$.
The selection strategy includes:
\begin{description}
 \item[Calculation with an integer step.] Using the explicit scheme, we determine the solution $v_{n+1}$
  at the time level $t_{n+1}$ via the step $\tau_{n+1}$;
 \item[Calculation with a half-integer step.] Using the explicit scheme, we calculate the solution $w_{n+1}$
  at the time level $t_{n+1}$ employing the step $0.5\tau_{n+1}$;
 \item[Step selecting.] The step $\tau_{n+1}$ is evaluated through the closeness between
  $v_{n+1}$ and $w_{n+1}$.
\end{description}  

For $v_{n+1}$, we have (\ref{3.2}), and $w_{n+1}$ is determined as follows:
\begin{equation}\label{4.1}
  \frac{w_{n+1/2} - y_{n}}{0.5 \tau_{n+1}} + A_{n} y_{n} = f_{n} .
\end{equation} 
\begin{equation}\label{4.2}
  \frac{w_{n+1} - w_{n+1/2}}{0.5 \tau_{n+1}} + A_{n+1/2} w_{n+1/2} = f_{n+1/2} .
\end{equation} 
Eliminating $w_{n+1/2}$ from (\ref{4.1}), (\ref{4.2}), we get
\[
\begin{split}
 w_{n+1} & = y_{n} - \frac{ \tau_{n+1}}{2} (A_{n+1/2} + A_n) y_n \\
 & +
 \frac{ \tau_{n+1}}{2} (f_{n+1/2} + f_n)  -
 \frac{ \tau_{n+1}^2}{4} A_{n+1/2}(f_n - A_n y_n) .
\end{split}
\]
Because of this, we have 
\begin{equation}\label{4.3}
\begin{split}
 v_{n+1} - w_{n+1} & =  \frac{ \tau_{n+1}}{2} (A_{n+1/2} - A_n) y_n \\
 & -  
 \frac{ \tau_{n+1}}{2} (f_{n+1/2} - f_n)  +
 \frac{ \tau_{n+1}^2}{4} A_{n+1/2}(f_n - A_n y_n) .
\end{split}
\end{equation} 
In view of the above notation (\ref{3.5}), we employ the approximate expressions:
\[
 \frac{ \tau_{n+1}}{2} (f_{n+1/2} - f_n) \approx \frac{\tau_{n+1}^2}{4 \gamma \tau_{n}} (\widetilde{f}_{n+1} - f_{n}) ,
\] 
\[
 \frac{ \tau_{n+1}}{2} (A_{n+1/2} - A_n) y_n \approx  \frac{\tau_{n+1}^2}{4 \gamma \tau_{n}} (\widetilde{A}_{n+1} - A_{n}) y_{n} . 
\]
By (\ref{2.4}), we have
\[
 \frac{ \tau_{n+1}^2}{4} A_{n+1/2}(f_n - A_n y_n)  \approx  \frac{\tau_{n+1}^2}{4 \tau_{n}}  \widetilde{A}_{n+1} (y_{n}-y_{n-1}) .
\]
Thus, we arrive at
\begin{equation}\label{4.4}
  v_{n+1} - w_{n+1} \approx \frac{\tau_{n+1}^2}{4 \tau_{n}} \left (  
  \frac{1}{\gamma} (\widetilde{A}_{n+1} - A_{n}) y_{n} - \frac{1}{\gamma} (\widetilde{f}_{n+1} - f_n) +
 \widetilde{A}_{n+1} (y_{n}-y_{n-1}) \right ) .
\end{equation}
The right-hand side of (\ref{4.4}) coincides with the right-hade side of (\ref{3.6}) with an accuracy of a factor. Similarly to (\ref{3.11}), we can formulate the rule for selecting the time step:
\begin{equation}\label{4.5}
  \tau_{n+1} = \max\left \{\tau_0, \min\{\gamma, \widetilde{\gamma}_{n+1}\} \tau_{n} \right \} ,
  \quad \widetilde{\gamma}_{n+1} = 4 \gamma_{n+1} .
\end{equation}
In fact, we have come to the same rule for the estimation of the time step --- the factor 4 has not
any essential matter.

\section{Numerical experiments}

To demonstrate the performance of the proposed algorithm (\ref{3.5}), (\ref{3.9}) for
selecting a time step based on the implicit scheme for solving the problem (\ref{2.1}), (\ref{2.2}), let us consider the boundary value problem for a one-dimensional parabolic equation.
Let $u(x,t)$  satisfies the equation
\begin{equation}\label{5.1}
  \frac{\partial u}{\partial t} - \frac{\partial^2 u}{\partial x^2} + p(t) u = f(t),
  \quad 0 < x < 1,
  \quad 0 < t \leq  T ,  
\end{equation}
as well as the boundary and the initial conditions: 
\begin{equation}\label{5.2}
  u(0, t) = 0,
  \quad u(1,t) = 0,
  \quad 0 < t \leq  T ,    
\end{equation} 
\begin{equation}\label{5.3}
  u(x,0) = u_0(x),
  \quad  \quad 0 <  x <  1 .
\end{equation} 

To solve approximately the problem (\ref{5.1})--(\ref{5.3}),
we apply finite difference discretization in space. Let us introduce 
a uniform grid with a step $h$:
\[
  \bar{\omega} = \{ x \ | \ x = ih, \quad i = 0, 1, ..., M, \quad Mh = 1 \},  
\] 
and $\omega$ is the set of interior grid points, whereas  $\partial \omega$ is the set of boundary points
($\bar{\omega } = \omega \cup \partial \omega$).
On the set of grid functions such that $u(x) = 0, \ x \notin \omega$, we introduce
a Hilbert space $H$, where the inner product and the norm are defined as:
\[
 (u,v) = \sum_{x \in \omega} u(x) v(x) h,
 \quad \|u\| = (u,u)^{1/2} . 
\]
The grid operator $A(t)$ is written as follows:
\[
 A u = - \frac{1}{h^2} (u(x+h) - 2 u(x) + u(x-h)) + p(t) u(x),
 \quad x \in \omega . 
\]
The operator $A(t)$ is self-adjoint, and if $p(t) \geq 0$, then it is positive definite in $H$.
Thus, after the spatial discretization of (\ref{5.1})--(\ref{5.3}), we arrive to
the problem (\ref{2.1}), (\ref{2.2}).

As a test problem, we consider the problem (\ref{5.1})--(\ref{5.3})
with $T = 0.1$ and the discontinuous coefficient $p(t)$  and  the discontinuous $f(t)$
defined as follows:
\[
  p(t) = \left \{
  \begin{split}
  100 t, & \quad 0 < t \leq 0.075, \\
  0, & \quad 0.075 < t \leq 0.1 , \\   
  \end{split}  
  \right . 
\]  
\[
  f(t) = \left \{
  \begin{split}
  0, & \quad 0 < t \leq 0.05, \\
  10 e^{-(t-0.05)}, & \quad 0.05 < t \leq 0.1 . \\   
  \end{split}  
  \right . 
\]
The problem is solved on the grid with $M = 100$, the calculations are performed using the sufficiently small time step $\tau_1 = \tau_0 = 1 \cdot 10^{-6}$.

\begin{figure}[htp]
  \begin{center}
    \includegraphics[scale = 0.5] {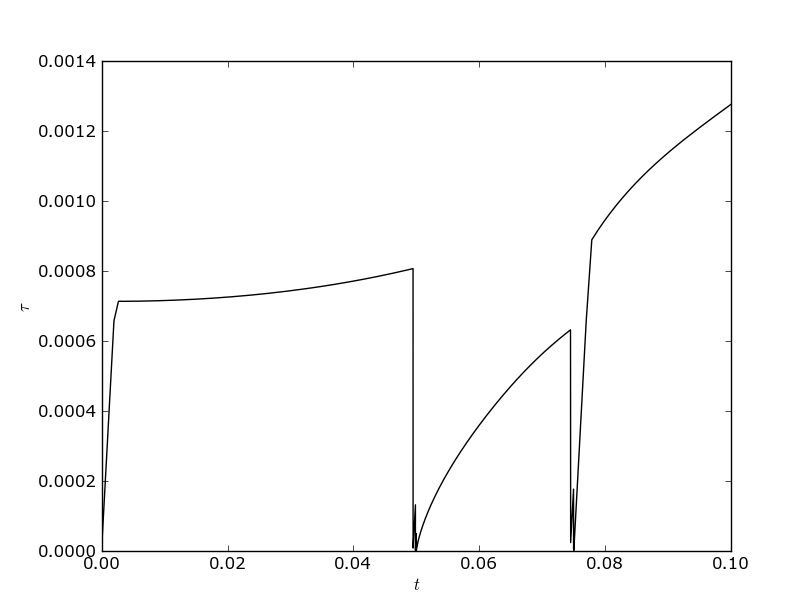}
	\caption{The time step: $\delta = 0.1, \ N = 268$ }
	\label{fig:1}
  \end{center}
\end{figure}

First, we Consider the case, where the initial condition  (\ref{5.3}) is taken in the following form:
\[
  u_0(x) = \sin(\pi x),
  \quad  \quad 0 <  x <  1 .
\]
If we specify the error level $\delta = 0.1$ and the parameter $\gamma =1.5$, then the time step 
produced by the algorithm (\ref{3.7}), (\ref{3.11}) has the form depicted in Fig.~\ref{fig:1}. The total number 
of time steps is $N = 268$.

\begin{figure}[htp]
  \begin{center}
    \includegraphics[scale = 0.5] {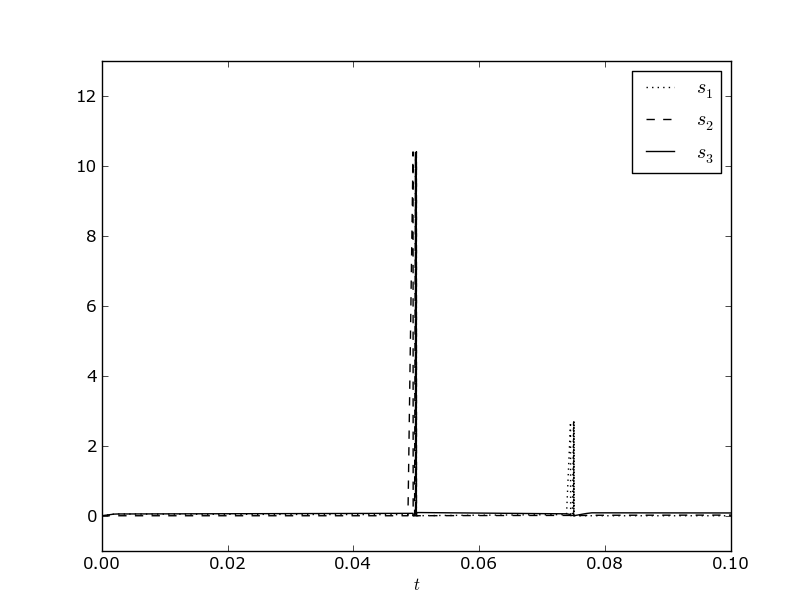}
	\caption{Terms of the correcting coefficient}
	\label{fig:2}
  \end{center}
\end{figure}

In this figure, we observe essential changes in the value the time step at $t=0.05$ and $t=0.075$, i.e.,
at the time moments corresponding to discontinuities in the right-hand side and the coefficient of the equation.  In accordance with the  rule (\ref{3.6}), the time step increases at the initial time stage.
Let us decompose the correcting coefficient $\gamma_{n+1}$ into three terms:
\[
 s_1 = \frac{1}{\gamma} \| (\widetilde{A}_{n+1} - A_{n}) y_{n} \|,
 \quad s_2 = \frac{1}{\gamma} \|\widetilde{f}_{n+1} - f_n \| ,
 \quad s_3 = \| \widetilde{A}_{n+1} (y_{n}-y_{n-1}) \| .
\]
Figure~\ref{fig:2} demonstrates their influence.

\begin{figure}[htp]
  \begin{center}
    \includegraphics[scale = 0.5] {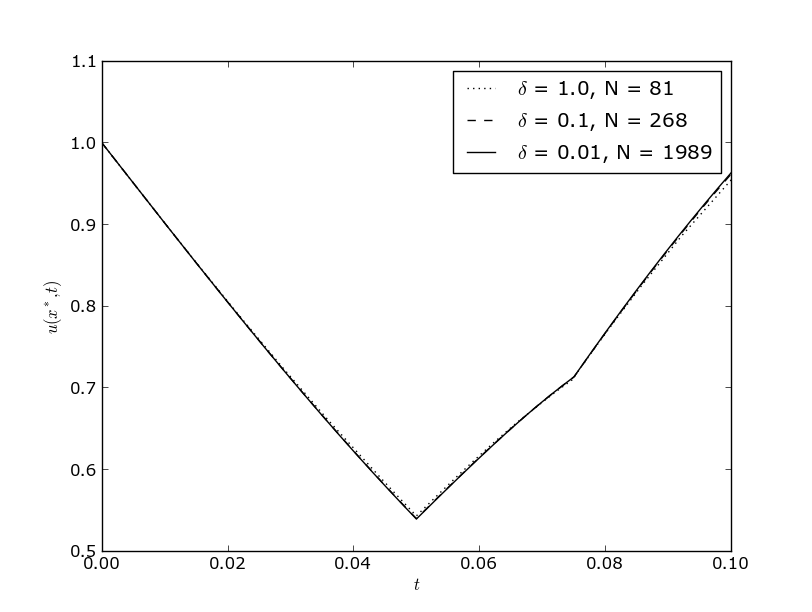}
	\caption{Approximate solution for various $\delta$}
	\label{fig:3}
  \end{center}
\end{figure}

\begin{figure}[htp]
  \begin{center}
    \includegraphics[scale = 0.5] {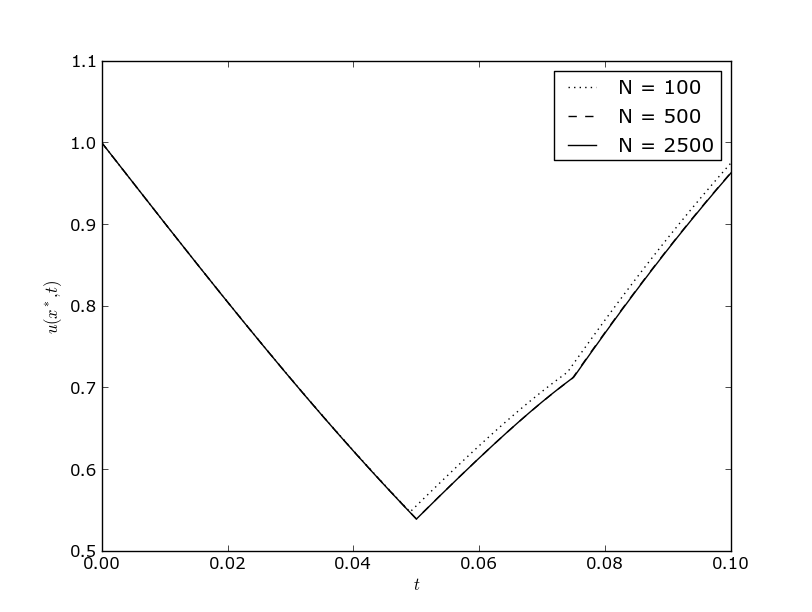}
	\caption{Approximate solution for various uniform time grids}
	\label{fig:3a}
  \end{center}
\end{figure}

The influence of the reducing error level $\delta$ on the convergence of the approximate solution
is shown in Fig.~\ref{fig:3}. The approximate solution at the point $x^* = 0.5$ is depicted in this figure.
For comparison, Figure\ref{fig:3a} presents similar data that were obtained using the uniform grids in time.

\begin{figure}[htp]
  \begin{center}
    \includegraphics[scale = 0.5] {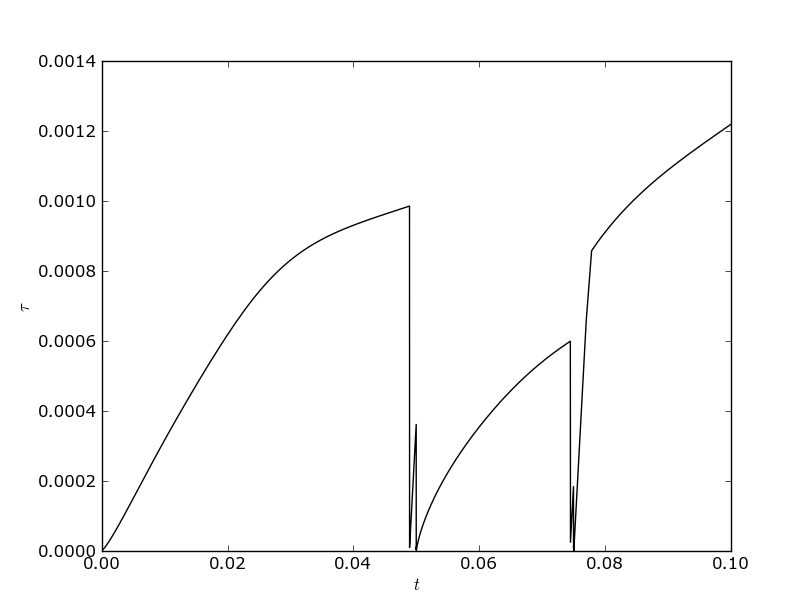}
	\caption{Time step: $\delta = 0.1, \ N = 569$ }
	\label{fig:4}
  \end{center}
\end{figure}

\begin{figure}[htp]
  \begin{center}
    \includegraphics[scale = 0.5] {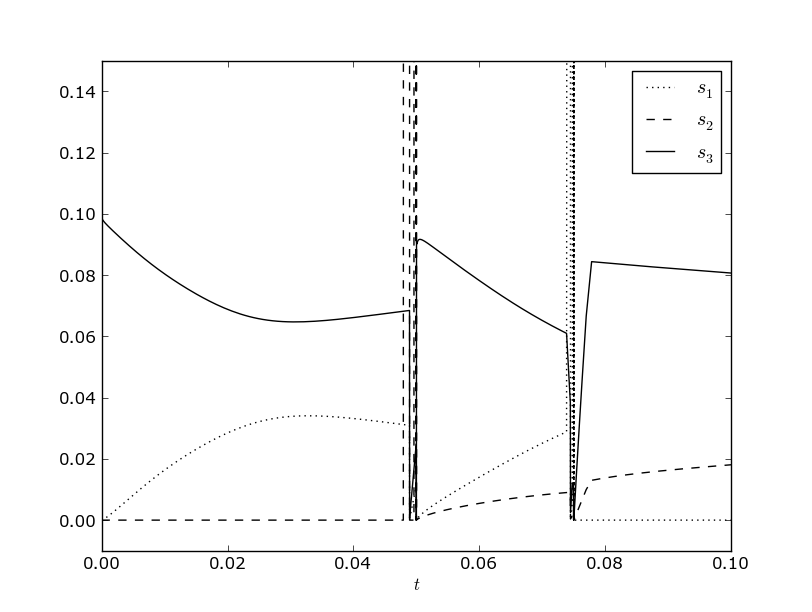}
	\caption{Correcting terms}
	\label{fig:5}
  \end{center}
\end{figure}

Special attention should be given to the influence of the initial conditions.
A typical situation is the presence of a boundary layer and this requires
to use small steps at the initial  time stage. For example, the behavior
of the time step for our model problem  with initial conditions
\[
  u_0(x) = \left \{
  \begin{split}
  2x, & \quad 0 < x \leq 0.5, \\
  1-2(x-0.5), & \quad 0.5 < x < 1 , \\   
  \end{split}  
  \right . 
\]
is shown in Fig.~\ref{fig:4}.
Compared with Fig.~\ref{fig:1} (smooth initial conditions),
the initial time stage is calculated with essentially smaller time steps
and the total number of steps is increased by more than a factor of 2.
In the region outside the neighbourhood of discontinuities of the coefficients and the right-hand side, the time step is controlled first of all by the term $s_3$ (see Fig.~\ref{fig:5}).

\begin{figure}[htp]
  \begin{center}
    \includegraphics[scale = 0.5] {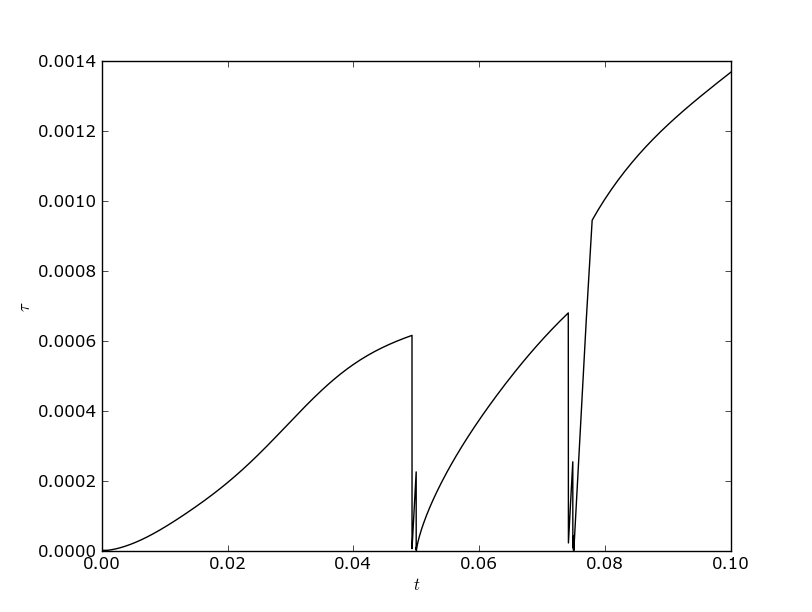}
	\caption{Time step: $u_0(x) = 1$}
	\label{fig:6}
  \end{center}
\end{figure}

A more difficult situation for the numerical solution is connected with inconsistent initial and boundary conditions. Let we have
\[
  u_0(x) = 1,
  \quad 0 < x < 1 . 
\]
The selection of the time step for this case is shown in Fig.~\ref{fig:6}.
Up to $t = 0.000856$ the calculation is carried out with the minimum time step $\tau_0 = 1 \cdot 10^{-6}$.
That is why the total number of time steps is 2183.

\end{document}